\documentstyle{amsppt}
\hoffset =.1in
\NoBlackBoxes
\magnification=\magstep 1
\nologo

\baselineskip18pt
\pageheight{25truecm}
\pagewidth{15truecm}
\vcorrection{-1truecm}

\font\erm=cmr8

\nopagenumbers

\define\Proj{\operatorname{Proj}}

\define\lin{\operatorname{lin}}


\topmatter
\title  The Amemiya-Ando conjecture falls
\endtitle
\author Adam Paszkiewicz
\endauthor

\abstract In any infinite dimensional Hilbert space $H$, a sequence $P_n\dots P_1x$ diverges in norm for some $x\in H$ and orthogonal projections $P_n\in\{Q_1,\dots, Q_5\}.$
\endabstract
\endtopmatter

\footnote""{{\it 2010} AMS subject classification: 47H09.}
\footnote""{Research supported by  grant N N201 605840.}

We prove the following theorem.

\proclaim{Theorem 0} In any infinite dimensional Hilbert space $H$, a sequence $P_n\dots P_1x$ diverges in norm for some $x\in H$ and orthogonal projections $P_n\in\{Q_1,\dots, Q_5\}.$

\endproclaim

In [1], a weak convergence of $P_n\dots P_1 x$ is proved for $P_n\in\{Q_1,\dots , Q_r\}$ with any $r$ "while strong convergence is still unsettled".  The norm convergence is known for $r =2$, for $r =3$ or $4$ the problem seems to be open. It is also known that for some projections $Q_1, Q_2$ on closed convex cones the   divergence  in norm of $(Q_1 Q_2)^n x$ can be obtained for any infinite dimensional $H$, and some $x\in H$ (see [2], [3], [4], [5]).

By capital letters $E, F, P, Q,\dots$ we shall always denote orthogonal projections in $H$. By small letters $e, f, p, q$ we shall always denote some vectors in $H$, of norm $1$; $e'$ is always orthogonal to $e$ ($f'\bot f$ and so on). Then $\hat e\cdot= \langle \cdot, e\rangle e$ is a one-dimensional projection and always $\hat e'\bot \hat e$, then $\hat f\leq \hat e +\hat e'$ means that $f\in \lin (e, e').$

The proof of Theorem 0 is unexpectedly simple. The proof of the following lemma is the most complicated part of the construction. 

\proclaim{Lemma 1} For any one-dimensional projections $\hat f_0,\dots, \hat f_r\leq \hat e +\hat e'$, and for $\epsilon >0$ there exist projections $P_0\leq \dots \leq P_r$ and numbers $n(0)\leq \dots \leq n(s)$ such that for $E =\hat e +\hat e'$
$$||(E P_s E)^{n(s)} - \hat f_s|| <\epsilon,\qquad 0 \leq s \leq r, \tag1$$
and one can require that $\dim(P_s - P_{s -1}) =1$ for $1\leq s\leq r.$

\endproclaim

\demo{Proof} Let us denote by $f'_s$ a vector (of norm $1$) satisfying  $\hat f_s +\hat f'_s = E$. For fixed mutually orthogonal vectors $p_0,\dots, p_r$ orthogonal to $E$, and for $\delta_0,\dots, \delta_r >0$ (defined later), we denote by $U_s$ a unitary operator
$$U_s f'_s = f'_s \cos \delta_s + p_s \sin \delta_s,$$
$$U_s p_s = - f'_s \sin \delta_s + p_s \cos \delta_s,$$
$$U_s x = x \quad\text{for } \ x\bot f'_s,  p_s .$$

Then for any projection $P\geq E$, orthogonal to $p_s$ for some $0\leq s \leq r$, and for $P_s =P - \hat f'_s$ we have $P_s E = \hat f_s$ and $U_s P U_s^* = \widehat{U_s f'_s} + P_s$ and thus
$$(E U_s P U_s^* E)^n =\hat f_s + (\hat f'_s \widehat{ U f'_s} \hat f'_s)^n =\hat s_s + \hat f'_s\cos^{2n} \delta_s\quad\text{for } \  n\geq 1.$$
In particular for any $0 < \eta <\frac{\pi}{2}$ there exists $n\geq 1$ satisfying
$$||(E U_s P U_s^* E)^n - \hat f_s|| <\eta \tag2$$
and $U_s P U_s^* \leq P + \hat p_s$ if only $P \geq E$, $P\bot \hat p_s$, $0\leq s\leq r$. The required projections are given by formula
$$P_s = U_r\dots U_s (E + \hat p_0 +\dots + \hat p_{s-1}) U_s^* \dots U_r^*,\qquad 0\leq s \leq r.$$
In fact,  we have $P_{s+1}\geq P_s$, $\dim(P_{s+1} - P_s) =1$, which follows from relations $E +\hat p_0 +\dots +\hat p_s \geq U_s(E +\hat p_0 +\dots + \hat p_{s-1})U_s^*$, $\dim(E +\hat p_0 +\dots +\hat p_s - U_s(E + \hat p_0 +\dots +p_{s-1})U_s^*)=1$, for any $0\leq \delta_s < \frac{\pi}{2}$, $0\leq s \leq r.$

The estimate (1) can be obtained for $\delta_s$, $n(s)$ chosen as  follows. The number $0< \delta_0 < \frac{\pi}{2}$ is taken arbitrarily and $n(0)$ is so large that, by (2),
$$||(E U_0 E U_0^* E)^{n(0)} - \hat f_0|| < (1 -\frac{1}{2^1})\epsilon.$$

Let $\delta_t >0$, $n(t)\geq 1$ be chosen for $0\leq t \leq s-1$, and let
$$||(E U_{s-1}\dots U_t(E + \hat p_0 +\dots + \hat p_{t-1})U^*_t \dots U_{s-1} E)^{n(t)} -\hat f_t|| < (1 - \frac{1}{2^s})\epsilon,\quad 0\leq t \leq s -1,$$
for a fixed $1 \leq s \leq r$. Then for sufficiently small $\delta_s > 0$ we have
$$\multline||(E U_s \dots U_{t}(E +\hat p_0 +\dots +\hat p_{t-1}) U^*_{t} \dots U^*_s E)^{n(t)} -\hat f_t|| \\f< (1 - \frac{1}{2^{s+1}})\epsilon,\quad 0\leq t \leq s -1,\endmultline\tag3$$
but for such $\delta_s$, there exists so large $n(s) \geq 1$ that, by (2),
$$||(E U_s(E +\hat p_0 +\dots + \hat p_{s-1})U_s^* E)^{n(s)} -\hat f_s|| < (1 - \frac{1}{2^{s+1}})\epsilon.$$
We have proved (3) for any $0\leq t\leq s$ and for $1\leq s\leq r.$

In particular, putting $s =r$ we have
$$||(E P_t E)^{n(t)} -\hat f_t|| < (1 -\frac{1}{2^{r+1}})\epsilon,\qquad 0\leq t \leq r.$$
\vskip-10pt
\hfill $\square$
\enddemo

\baselineskip18pt

\proclaim{Lemma 2} For any projections $P_0\leq \dots \leq P_r$ $\dim(P_s - P_{s-1}) =1$, $1\leq s \leq r$, and for $\epsilon >0$, there exist $Q\in \Proj H$ and numbers $m(0),\dots, m(r) \geq 1$ satisfying
$$||(P_r Q P_r)^{m(s)} - P_s|| < \epsilon, \qquad 1\leq s\leq r. \tag4$$

\endproclaim

\demo{Proof} Let $P_s = P_0 +\hat p_1 +\dots + \hat p_s$, and $P_0, \hat p_1,\dots, \hat p_r, \hat p'_1,\dots, \hat p'_r$ be mutually orthogonal. For some $\delta_1,\dots, \delta_r >0$, let $q_s =p_s \cos \delta_s +p'_s \sin \delta_s$ and let $Q = P_0 +\hat q_1 + \dots + \hat q_r$. Now, to obtain (4), we  can proceed as follows. We take $m(r) =1$ and so small $\delta_r >0$, that $||(\hat p_r \hat q_r \hat p_r)^{m(r)} - \hat p_r|| <\epsilon$. Then we proceed downward. Assuming that for some $1 < s \leq r$ and $m(s) >\dots > m(r)$, $0 < \delta_s, \dots, < \delta_r$ we have
$$\aligned &||(\hat p_s \hat q_s \hat p_s)^{m(t)} -\hat p_s|| <\epsilon,\qquad s\leq t \leq r,\\
&||(\hat p_t \hat q_t \hat p_t)^{m(s)}|| < \epsilon, \qquad \qquad s+1 \leq t\leq r,\endaligned\tag5$$
it is possible to find such $m(s-1) > m(s)$ that $||(\hat p_t \hat q_t \hat p_t)^{m(s-1)}|| < \epsilon$ for $t =s$ (and thus for any $s \leq t\leq r$).  Then there exists so small $\delta_{s-1}$, $0\leq \delta_{s-1} < \delta_s$, that
$$||(\hat p_{s-1} \hat q_{s-1} \hat p_{s-1})^{m(t)} - \hat p_{s-1}|| <\epsilon \qquad \text{for } \ t = s-1$$
(and thus for any $s-1 \leq t \leq r$).

Now (5) is proved for any $1\leq s\leq r$, and this implies (4) if only $m(0) > m(1)$ satisfies $||(\hat p_1 \hat q_1\hat p_1)^{m(0)}|| <\epsilon.$

\remark{Remark} The assumption $\dim (P_s - P_{s-1}) =1$ makes the proof of Lemma 2 more elementary, but is not essential.

\endremark

\proclaim{Corollary 3} For any one-dimensional projections $\hat f_0,\dots \hat f_r \leq \hat e +\hat e'$ and for $\epsilon >0$ there exist projections $P, Q$, and numbers $m(s), n(s) \geq 1$, $0\leq s\leq r$ satisfying
$$||(E(PQP)^{m(s)} E)^{n(s)} -\hat f_s|| <\epsilon,$$
for $E =\hat e +\hat e'.$

\endproclaim

\demo{Proof} By Lemma 1, the estimate (1) can be satisfied  for some $P_0 \leq \dots \leq P_r$, $\dim(P_s - P_{s-1}) =1$, and for $\epsilon /2$ instead of $\epsilon$. Then we have
$$||(E(P_r Q P_r)^{m(s)} E)^{n(s)} - \hat f_s|| < \epsilon,\quad 0\leq s\leq r,$$
if only (4) is satisfied for sufficiently small $\epsilon_1 >0$ instead of $\epsilon$. Lemma 2 finishes the proof.

\enddemo

\proclaim{Lemma 4} For any vectors $e\bot e'$ and $\epsilon >0$, there exist projections $P, Q$ and a monomial $A(P, Q, E) = P_k \dots P_1$, with $P_l \in \{P, Q, E\}$, $1\leq l\leq k$,  $E = \hat e+ \hat e'$, such that
$$\langle A(P, Q, E)e, e'\rangle >1 -\epsilon.$$

\endproclaim

\demo{Proof} For $f_s = e\cos \frac{\pi s}{2 r} + e' \sin \frac{\pi s}{2 r}$, $0\leq s\leq r$, an elementary calculation gives
$$\langle \hat f_r\dots \hat f_0 e, e' \rangle = (\cos \frac{\pi}{2r})^r > 1 - \frac{\epsilon}{2} \quad\text{for some } \ r\geq 1.$$
Thus Corollary 3 (with sufficiently small $\epsilon_1 >0$ instead of $\epsilon$) gives
$$\langle (E(P Q P)^{m(r)} E)^{n(r)} \dots (E(PQP)^{m(0)} E)^{n(0)} e, e'\rangle \geq 1 -\epsilon.$$

\enddemo

\flushpar{\bf 5. Proof of Theorem 0}. Let us take infinitely dimensional projections $F_1, F_2,\dots$ satisfying, for some orthonormal ($e_i$,  $i\geq 1$),
$$\hat e_1 \leq F_1, \quad F_i F_{i+1} =\hat e_{i+1},\qquad i\geq 1,$$
$$F_i\bot F_j\qquad\text{for } \ |i - j| \geq 2$$
and denote additionally  $e'_i = e_{i+1}$, $E_i =\hat e_i +\hat e'_i \leq F_i$, $i\geq 1.$

For fixed $0 <\epsilon_i <1$, $\prod (1-\epsilon_i) >0$, we use Lemma 4 in subspaces $H_i = F_i(H)$ (instead of $H$) obtaining projections $P_i, Q_i \leq F_i$ and monomials $A_i(P_i, Q_i, E_i)$ satisfying
$$\eta_i:= \langle A_i(P_i, Q_i, E_i) e_i, e'_i\rangle  > 1 - \epsilon_i, \quad i\geq 1.\tag6$$

Let $P =\sum P_{2k}$, $Q = \sum Q_{2k}$, $E = \hat e_1 +\sum E_{2k}$, for $k\geq 1$, then $E =\sum \hat e_i$, $E F_{2k} = E_{2k}$ and for any monomial $A(P, Q, E)$ we have
$$A(P, Q, E)F_{2k}(H)\subset F_{2k}(H),$$
$$A(P, Q, E)F_{2k-1} = A(P, Q, E)(\hat e_{2k-1} +\hat e_{2k}) = A(P, Q, E)\hat e_{2k-1} + A(P_{2k}, Q_{2k}, E_{2k})\hat e_{2k}.$$
In consequence
$$\aligned  & x\in F_1\vee \dots \vee F_{2k-1}(H), \quad \langle x, e_{2k}\rangle =\eta\\
& \text{imply} \ \ A(P, Q, E)x  = y + \eta A(P_{2k}, Q_{2k}, E_{2k})e_{2k} \\ 
& \text{for some} \quad  y \in F_1 \vee \dots \vee F_{2k-2}(H),\endaligned$$
obviously $y =0$ for $k =1$.

Analogically, for $R =\sum P_{2k-1}$, $S =\sum Q_{2k-1}$, $E =\sum E_{2k-1}$, $k\geq 1$, (thus as before $E =\sum \hat e_i$) we have
$$\aligned  & x\in F_1\vee \dots \vee F_{2k}(H), \quad \langle x, e_{2k+1}\rangle =\eta\\
& \text{imply} \ \ A(R, S, E)x  = y + \eta A(R_{2k+1}, Q_{2k+1}, E_{2k+1})e_{2k+1}\\
& \text{for some} \quad  y \in F_1 \vee \dots \vee F_{2k-1}(H).\endaligned$$

Let us denote
$$x_0 =e_1, \quad x_1 = A_1(R, S, E)x_0$$
and, by induction for $k\geq 1$,
$$x_{2k} = A_{2k}(P, Q, E)x_{2k-1},\quad x_{2k+1} = A_{2k+1}(R, S, E)x_{2k}.$$
Then (6) implies
$$x_1\in F_1(H),\qquad \langle x_1, e_2\rangle =\eta_1,$$
and for all $k\geq 1$
$$x_{2k-1}\in F_1\vee \dots\vee F_{2k-1}(H),\qquad \langle x_{2k-1}, e_{2k}\rangle =\eta_1\dots \eta_{2k-1},$$
$$x_{2k}\in F_1\vee \dots \vee F_{2k}(H),\qquad \langle x_{2k}, e_{2k+1}\rangle = \eta_1\dots \eta_{2k}.$$
As $\prod \eta_n > \prod (1 -\epsilon_i) >0$, the sequence
$$A_{2k}(P, Q, E) A_{2k-1}(R, S, E) \dots A_2(P, Q, E) A_1(R, S, E)e_1$$
diverges in norm. Thus $\{Q_1, \dots, Q_5\} = \{P, Q, R, S, E\}$ can be taken.

\bigskip

\widestnumber\key{AAA}

\flushpar{\bf References}
\medskip

\ref\key{1}
\by I. Amemiya, T. Ando
\paper Convergence of random products of contractions in Hilbert space
\jour Acta Sci. Math. (Szeged)
\vol 26
\yr 1965
\pages 239-244
\endref

\ref\key{2}
\by H.H. Bauschke, E. Mateo\v skov\'a, S. Reich
\paper Projections and proximal point methods: Convergence results and counterexamples
\jour Nonlinear Anal.
\vol 56
\yr 2004
\pages 715-738
\endref

\ref\key{3}
\by H.S. Hundal
\paper An alternating projection that does not converge in norm
\jour Nonlinear Anal.
\vol 57
\yr 2004
\pages 35-61
\endref

\ref\key{4}
\by E. Kopeck\'a
\paper Spokes,  mirrors and alternating projections
\jour Nonlinear Anal.
\vol 68
\yr 2008
\pages  1759-1764
\endref

\ref\key{5}
\by E. Matou\v skov\'a, S. Reich
\paper The Hundal example revisited
\jour J. Nonlinear Convex Anal.
\vol 4
\yr 2003
\pages 411-427
\endref

\bigskip
\erm
\baselineskip12pt
\flushpar Faculty of Mathematics and Computer Science
\flushpar \L \'od\'z \,University
\flushpar Banacha 22, 90-238 \L \'od\'z, Poland
\flushpar e-mail: ktpis\@math.uni.lodz.pl

\bye